\documentclass[12pt]{article}
\usepackage{amssymb}
\usepackage{booktabs}
\def\inh{\vskip 0.075truein \noindent\hangindent=12 pt \hangafter=1}

\usepackage{amsthm}\theoremstyle{remark}

\newcommand{\bte}{\begin{quote}\begin{theorem}}
\newcommand{\ete}[1]{\label{#1}\end{theorem}\end{quote}}
\newcommand{\bcom}{\begin{quote}\end{quote}}
\newcommand{\bex}{\begin{quote}\begin{example}}
\newcommand{\eex}[1]{\label{#1}\end{example}\end{quote}}
\newcommand{\bcon}{\begin{quote}\begin{conclusion}}
\newcommand{\econ}[1]{\label{#1}\end{conclusion}\end{quote}}
\newcommand{\bdefi}{\begin{quote}\begin{definition}}
\newcommand{\edefi}[1]{\label{#1}\end{definition}\end{quote}}

\newcommand{\blem}{\begin{quote}\begin{lemma}}
\newcommand{\elem}[1]{\label{#1}\end{lemma}\end{quote}}

\newcommand{\bpr}{\begin{quote}\begin{problem}}
\newcommand{\epr}[1]{\label{#1}\end{problem}\end{quote}}

\newcommand{\f}{\frac}
\newcommand{\p}{\partial}
\newcommand{\n}{\nonumber \\}

\newcommand{\beq}{\begin{eqnarray}}
\newcommand{\eeq}[1]{\label{#1}\end{eqnarray}}
\newcommand\eq[1]{(\ref{#1})}
\newcommand{\bfi}{\begin{figure}[24]}
\newcommand{\efi}[1]{\caption{\label{#1}}\end{figure}}
\newcommand\fig[1]{Fig.~\ref{#1}}
\newcommand{\res}{respectively}
\newcommand\gl{\left}
\newcommand\gr{\right}
\newcommand{\bfm}[1]{\mbox{\boldmath ${#1}$}}

\newcommand{\Ga}{\alpha}

\newcommand{\Gf}{\phi}

\newcommand{\Gg}{\gamma}

\newcommand{\Gl}{\lambda}

\newcommand{\Gr}{\varrho}

\newcommand{\Go}{\omega}

\newcommand{\GD}{\Delta}

\newcommand{\GO}{\Omega}

\newcommand\D{\,\mathrm{d}}
\newcommand\I{\mathrm{i}}
\newcommand\E{\mathrm{e}}
\newcommand{\bexe}{\begin{quote}\begin{exercise}\inh}
\newcommand{\eexe}[1]{\label{#1}\end{exercise}\end{quote}}

\usepackage{graphics,graphicx}
\usepackage{color}
\begin{document}

\title{Resonant waves in elastic structured media: dynamic homogenisation versus Green's functions}

\author{{\small
Alexander B. Movchan${}^{a*}$, Leonid I. Slepyan${}^{b}$}\\
{\small ${}^a$ Department of Mathematical Sciences,
University of Liverpool, U.K.}\\
{\small ${}^b$ School of Mechanical Engineering, Tel Aviv
University, Israel}\\
{\small Institute of Mathematics and Physics,
Aberystwyth University, UK}}

\date{}

\maketitle
\begin{abstract}
We address an important issue of a dynamic homogenisation in vector elasticity for a doubly periodic  mass-spring elastic lattice. The notion of logarithmically growing resonant waves is used in a complete analysis of star-shaped wave forms induced by an oscillating point force. We note that the dispersion surfaces for Floquet-Bloch waves in an elastic lattice main contain critical points of the saddle type. Based on the local quadratic approximations of the frequency, as a function of wave vector components, we deduce properties of a transient asymptotic solution as the contribution of the point source to the wave form. In this way, we describe  local Green's functions  as localized wave forms corresponding to the resonant frequency. The peculiarity of the problem lies in the fact that, at the same resonant frequency, the Taylor quadratic approximations for different groups of the resonant points are different, and hence we deal with different local Green's functions. Thus, 
 there is no uniformly defined homogenisation procedure for a given resonant frequency. Instead, the  continuous approximation of the wave field can be obtained through the asymptotic analysis of the lattice Green's functions.

\end{abstract}

\vspace*{5mm}
 \noindent $^*$Email: abm@liverpool.ac.uk

\maketitle

\section{Introduction}

The subject of homogenisation is of great interest to physicists, engineers and mathematicians. The work in this area goes back more than 100 years, with the classical and elegant paper by Rayleigh (1892) being one of foundation stones in analysis of effective properties of periodic composite media. Mathematical theory of multi-scale homogenisation approximations has received a substantial attention, as comprehensively descibed  in the books by Bensoussan, Lions and Papanicoau (1978), Sanchez-Palencia (1980), Marchenko and Khruslov (2006), Bakhvalov and Panasenko (1984), and Zhikov, Zhukov et al (1994).

Conventional homogenisation in problems of wave propagation would normally apply to the case of long wave asymptotics, where a characteristic size of scatterers within a periodic structure is much smaller compared to the wavelength of the incident wave.

Resonant waves excited by a harmonic force in {\em uniform} square and triangular lattices were studied by Ayzenberg-Stepanenko  and Slepyan (2008) in the framework of a scalar problem. It was shown that the resonant waves spread mainly on some separate rays to form star-like configurations and hence show strong dynamic anisotropy. The underline structure dynamics was revealed and an asymptotic solution showing the localization was presented. Then the scalar problem for the {\em nonuniform}, periodic square lattice was examined by Craster et al (2010) (also see Craster et al (2009)), where the localization phenomena were also found for the scalar problem. The Craster et al (2009, 2010) papers also addressed the issue of a high-frequency homogenization in the neighbourhood of the resonant modes.

Recent publications by Milton et al. (2006), Milton and Nicorovici (2006), Norris (2008), Brun et al. (2009) on the dynamic response of metamaterials and invisibility cloaks in elastic and aciustic media raised interesting questions related to approximating of such systems by multi-scale composites, that may possess such unusual properties as chirality, negative refraction and negative inertia. The range of frequencies in such applications would be well outside the standard homogenisation range, and hence the new approach is required.

The paper by Colquitt et al. (2011) addresses the vector problems for elastic Floquet-Bloch waves in a beam-made triangular lattice. The analysis of the dispersion relations has revealed a strong dynamic anisotropy with a certain range of frequencies. It has also shown the effects of negative refraction for a certain class of structured interfaces.

The concept of high-frequency homogenisation for vector elasticity problems remains a challenge, as several critical points on a dispersion surface may exist for the same frequency. We address here these important issues for frequencies corresponding to resonant modes within the   mass-spring  triangular elastic lattice. Also, we focus on analysis of lattice Green's tensors for a two-dimensional vector problem of elasticity and show how this powerful approach compares to multi-scale homogenisation approximations. The analysis incorporates asymptotic approximations for the resonant waves. The directional localisation is then associated with saddle points on the dispersion surfaces, which are also  linked to the anisotropy in the dynamic regime.
According to the structure of the dispersion relations, there exist     several different wave forms corresponding to the same resonant frequency. Hence, different asymptotic solutions and homogenised equations may be derived for the same frequency. The subtlety is resolved by identifying the reference resonance modes and analysing the asymptotics of Green's tensors. We also identify asymptotic solutions for logarithmically growing resonance waves, which, in the case of a saddle point, also resemble the star waves.

\section{Governing equations for a forced elastic lattice}
Consider a regular triangular lattice consisting of point masses connected by  massless elastic links. The mass value, the bond length and stiffness are taken as the natural units. Thus, the lattice spacing along the bond line
is equal to unity, whereas the distance between the parallel bond lines is $\sqrt{3}/2$.
The lattice is subjected to a time-harmonic external force acting on a given mass. We would like to identify resonance modes, in particular those on the boundaries of the stop bands, and furthermore analyse homogenisation approximations corresponding to the resonant frequencies.

For the chosen geometry and physical parameters of the lattice,   in the long-wave/low-frequency
approximation, we obtain a homogeneous, isotropic, elastic body with the
density $\Gr=2/\sqrt{3}$, Poisson's ratio $\nu=1/3$ and the following velocities of the longitudinal,
shear and Rayleigh waves: $c_1 = \sqrt{9/8}$, $c_2 = \sqrt{3/8}$ and $c_R = \frac{1}{2}\sqrt{3-\sqrt{3}}$, respectively. The effective shear modulus for such a lattice in the static approximation is $\mu= \sqrt{3}/4$.

It will be shown that in the time-harmonic regime, as the frequency increases,  this lattice becomes far from being isotropic, and moreover, governing equations of the homogenised body may change its type from being elliptic to hyperbolic.

The displacement vector  is denoted by $\bfm{u}_{m,n}=(u_x,u_y)_{m,n}$, where the integers, $m,n$, define the mass position. In the Cartesian coordinates, $x,y,$ we have
 \beq x= m + \f{n}{2}\,,~~~y= \f{\sqrt{3}}{2}n\,;~~~(m,n)=0,\pm 1, \pm 2, ...\, .\eeq{mn}
 The displacement components satisfy the equations of motion as follows
 \beq \ddot{u}_{m,n}&& = Q_0-Q_3+\f{1}{2}(Q_1-Q_2-Q_4+Q_5) +P_{x;m,n}\,,\n
 \ddot{v}_{m,n} &&= \f{\sqrt{3}}{2}(Q_1+Q_2-Q_4-Q_5) +P_{y;m,n}\,,\eeq{1}
where $\bfm{P}_{m,n}=(P_{x;m,n},P_{y;m,n})$ are the external forces, applied to the $(m,n)$ mass,  whereas  $Q_j, j=0,1,..., 5,$ are the forces acting on the mass $(m,n)$ from the neighbouring masses, i.e.
 \beq Q_0&=&u_{m+1,n}-u_{m,n}\,,~~~Q_3=u_{m,n}-u_{m-1,n}\,,\n
 Q_1&=&\f{1}{2}(u_{m,n+1}-u_{m,n})+\f{\sqrt{3}}{2}(v_{m,n+1}-v_{m,n})\,,\n
 Q_2&=&-\f{1}{2}(u_{m-1,n+1}-u_{m,n})+\f{\sqrt{3}}{2}(v_{m-1,n+1}-v_{m,n})\,,\n
 Q_4&=&-\f{1}{2}(u_{m,n-1}-u_{m,n})-\f{\sqrt{3}}{2}(v_{m,n-1}-v_{m,n})\,,\n
 Q_5&=&\f{1}{2}(u_{m+1,n-1}-u_{m,n})-\f{\sqrt{3}}{2}(v_{m+1,n-1}-v_{m,n})\,.\eeq{2}

Assuming the time-harmonic vibrations with the amplitudes $\bfm{U}_{m,n}$ and radian frequency $\Go$, we have
$$\bfm{u}_{m,n}(t)=\bfm{U}_{m,n}\E^{\I\Go t}.$$
Correspondingly, the discrete Fourier transform gives
 \beq  \bfm{U}^{FF} (\bfm{k}) = \sum_{m,n}\bfm{U}_{m,n}\exp(\I\bfm{k}\cdot\bfm{x}(m,n))\,,\eeq{8}
where $\bfm{x}(m,n)$ is the position vector of the $(m,n)$-mass. The same transform applies to the amplitude of the external force acting on the masses.

The original amplitudes are then defined by the inverse transform, which in our particular case, is given in the form
 \beq \bfm{U}_{m,n} =  \f{\sqrt{3}}{16\pi^2}\int_{-2\pi/\sqrt{3}}^{2\pi/\sqrt{3}}\gl(\int_{-2\pi}^{2\pi} \bfm{U}^{FF} (\bfm{k})\E^{-\I(\bfm{k} \cdot \bfm{x}(m,n))}\D k_x\gr)\D k_y\,.\eeq{11}

\section{Elastic compliance and dispersion}

\begin{figure}[!ht]

\centering
\vspace*{0mm} \rotatebox{0}{\resizebox{!}{10cm}{%
\includegraphics [scale=0.25]{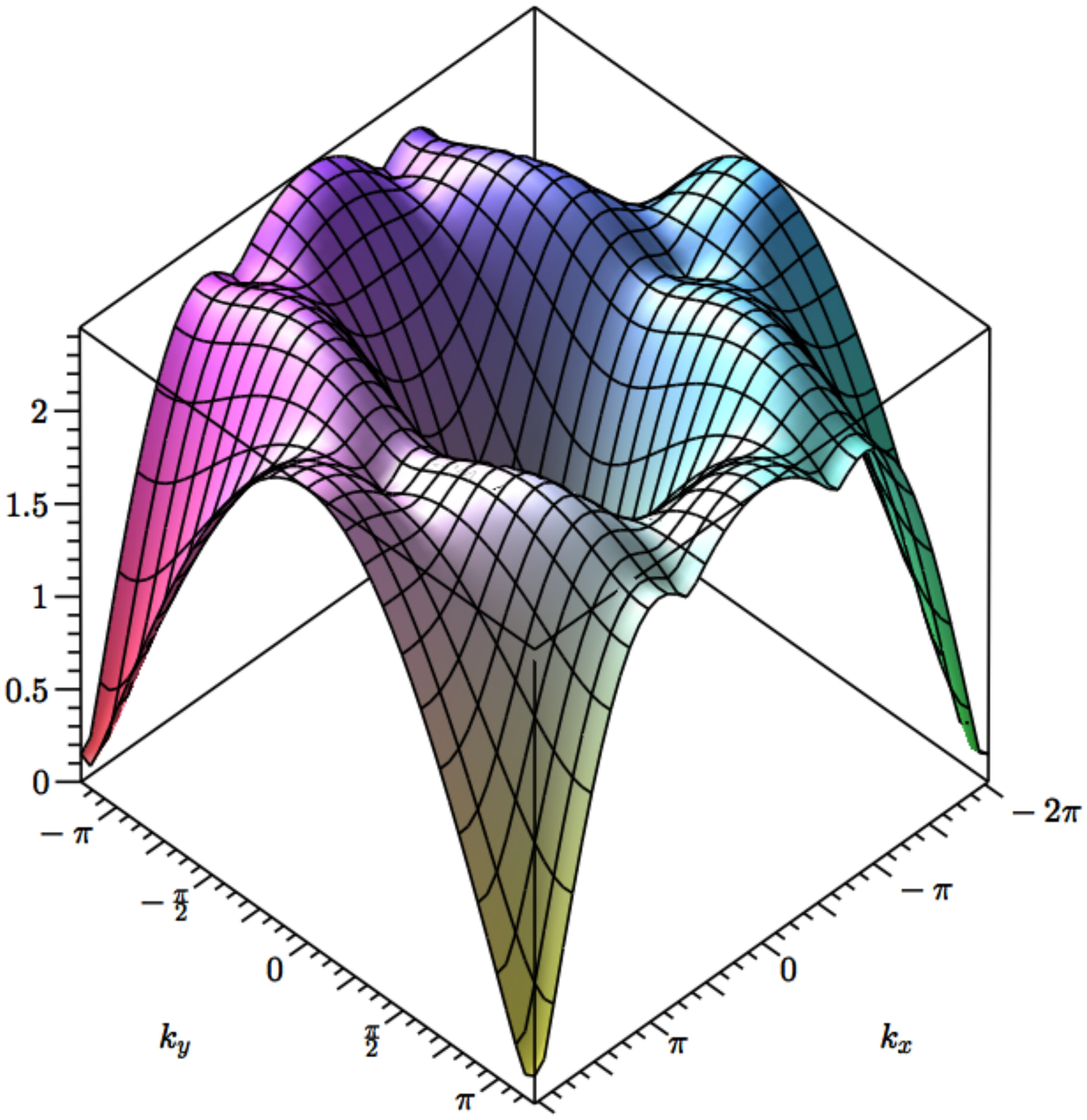}}}

\vspace{-5mm}
 \caption{Dispersion surfaces $\Go=\Go_1(k_x,k_y)$}
    \label{omega1}
\end{figure}

\begin{figure}[!ht]

\centering
\vspace*{-0mm} \rotatebox{0}{\resizebox{!}{10
cm}{%
\includegraphics [scale=0.25]{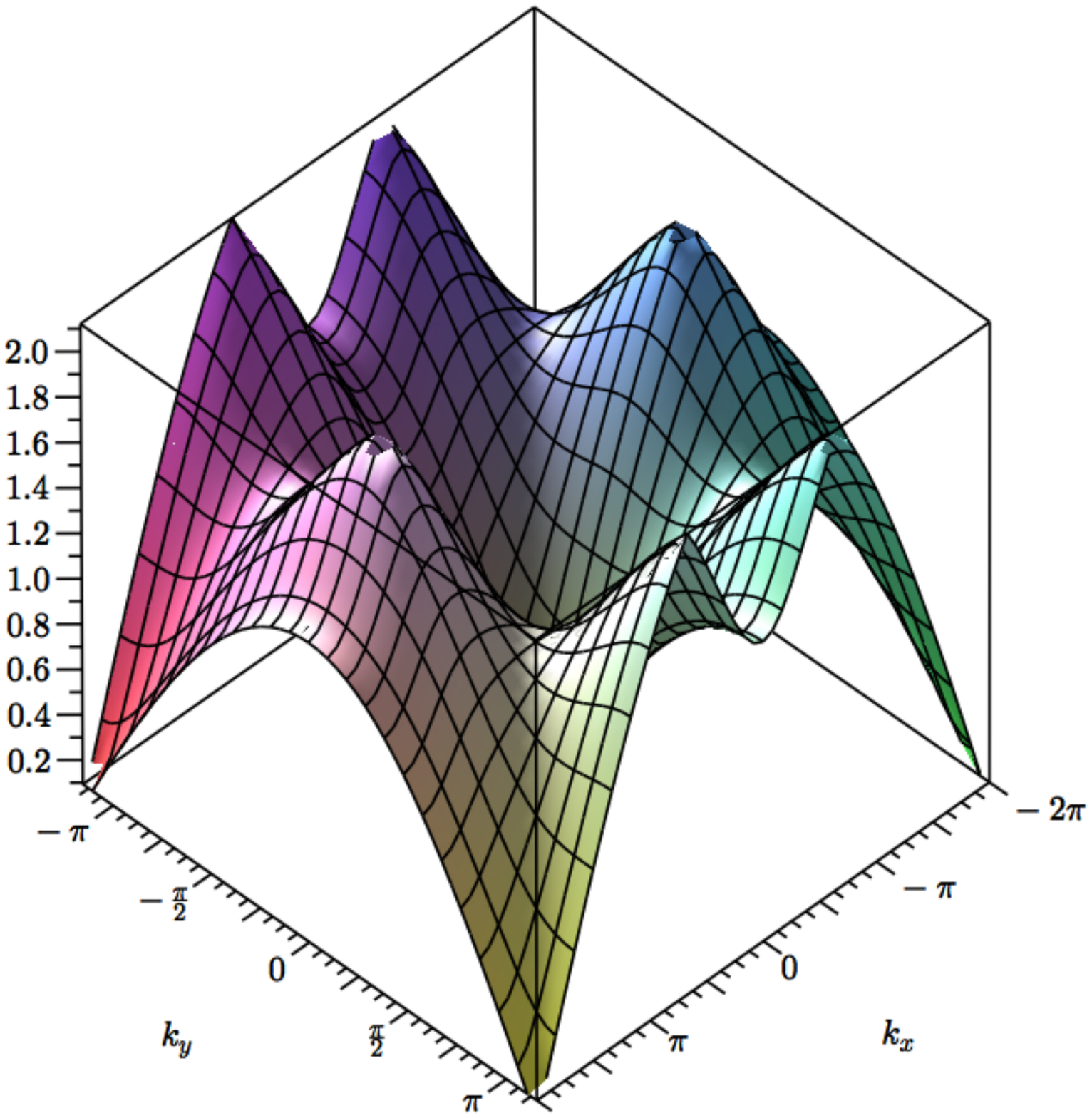}}}

\vspace{-0mm}
 \caption{Dispersion surface $\Go=\Go_2(k_x,k_y)$}
    \label{omega2}
\end{figure}

We assume that the external load is represented by a  time-harmonic point force with the amplitude vector $\bfm{P}_{0,0}=\bfm{P}=(P_x,P_y), |\bfm{P}|=1,$ acting on the central mass, $m=n=0$. With this in mind we find from \eq{1} that
 \beq \bfm{U}^{FF} = \bfm{A}_{\Go,\bfm{k}}\bfm{P}\,,\eeq{9}
 where $ \bfm{A}_{\Go,\bfm{k}}$ is the compliance symmetric tensor, given by
 \beq  \bfm{A}_{\Go,\bfm{k}} =\bfm{B}_{\Go,\bfm{k}}/\GD\,.
 \eeq{5a}
 Here
 \beq
  B_{xx}&=&-\GO+3\gl(1-\cos\gl(\f{k_x}{2}\gr)\cos\gl(\f{\sqrt{3}}{2}k_y\gr)\gr)\,,\n B_{xy}&=&B_{yx}= - \sqrt{3}\sin\gl(\f{k_x}{2}\gr)\sin\gl(\f{\sqrt{3}}{2}k_y\gr)\,,\n
 B_{yy} &=& -\GO+2(1-\cos k_x)+1-\cos\gl(\f{k_x}{2}\gr)\cos\gl(\f{\sqrt{3}}{2}k_y\gr)\,,   ~ \GO=\omega^2 , \eeq{10}
 and the function $\GD(\omega, k_x, k_y)$ 
 can be factorized in the form
 \beq \GD = (\GO -\GO_1)(\GO-\GO_2)\,,~~~\GO_{1,2} = \Go_{1,2}^2= F\pm\sqrt{S}\,.\eeq{5b}
 The quantities $F$ and $S$ are periodic functions of $k_x$ and $k_y$ defined as follows
 \beq
 F&=&1-\cos k_x+2\gl(1-\cos\gl(\f{k_x}{2}\gr)\cos\gl(\f{\sqrt{3}}{2}k_y\gr)\gr)\,,\n
 S&=& \gl(\cos k_x-\cos\gl(\f{k_x}{2}\gr)\cos\gl(\f{\sqrt{3}}{2}k_y\gr)\gr)^2 + 3\sin^2\gl(\f{k_x}{2}\gr)\sin^2\gl(\f{\sqrt{3}}{2}k_y\gr)\,.\eeq{5}
 The dispersion of waves in the elastic lattice is governed by the equation
 \beq
 \GD(\omega, k_x, k_y) =0.
 \eeq{5c}
Two dispersion surfaces, periodic in $k_x$ and $k_y$, are identified by the equations
$$\GO-\GO_{1,2} =0. $$
On the elementary cell of periodicity, these surfaces have common points, where $\GO_1=\GO_2= 0$, at the origin, $k_x=k_y=0$, and 
at the corner points
$k_x=\pm 2\pi, k_y=\pm 2\pi/\sqrt{3}$, where there is also $\GO_{1,2}= 0$.
The graphs of the dispersion surfaces $\Go=\Go_{1,2}(k_x,k_y)$ restricted  to the elementary cell of periodicity are presented in \fig{omega1} and \fig{omega2}.

\section{Resonant excitation. Asymptotics of Green's kernels}
Floquet-Bloch waves and their dispersion properties are important for evaluation of resonant frequencies within a periodic system. We refer to Ayzenberg-Stepanenko and  Slepyan  (2008), and Colquitt et al. (2011) who have analysed the resonant forms and star-shaped waves in lattices as well as features of dynamic anisotropy for elastic waves in an elastic triangular lattice.

\subsection{Critical points on the dispersion surfaces}

  The
  critical points on the dispersion surfaces  correspond to standing waves, and are identified by  the resonant frequencies, $\Go^*_i$, and components of the Bloch vector, as shown in  Table 1 below.

\begin{table}[htb]
\begin{tabular}{@{\qquad}c@{\qquad}c@{\qquad}c@{\qquad}c@{\qquad}c@{\qquad}}
\toprule
\multicolumn{1}{c}{} \\
$i$, ~   $\Go_i^*$  & $(k_{x;i1},k _{y;i1})$ &  $(k_{x;i2}, k_{y;i2})$ & $k_{x;i3}, k_{y;i3}$\\
\midrule
\multicolumn{1}{c}{
}\\
 $1$,~ $\sqrt{6}$ & $(0,  \pm 2\pi/\sqrt{3})$ & $(\pm 2\pi, 0)$ & $(\pm\pi, \pm 1.8137994)$ \\
 \midrule
 $2$,~   $2.25$ & $(\pm 2.8909370,\pm 2\pi/\sqrt{3})$ & $(\pm 3.3922483, 0)$ & $(\pm 1.6961242, \pm 2.9377732)$\\
  & & & $(\pm 4.5870611, \pm 0.6898255)$ \\
 \midrule
$3$,~  $\sqrt{2}$ & $(0, \pm 2\pi/\sqrt{3})$ & $(\pm 2\pi, 0)$ & $(\pm\pi, \pm\pi/\sqrt{3})$ \\
\bottomrule
\end{tabular}
\caption{\label{tab:error-measure-continuum}
Positions of critical points on the dispersion surfaces, as defined by components of the Bloch vector $(k_x, k_y)$ and the radian frequency $\Go^*$.}
\end{table}


We note that there are multiple critical points identified for the same frequency, which makes classical homogenisation  impossible.
 However, we intend to analyse the profile of the dispersion surface in the neighbourhoods of individual critical points and hence identify standing waves and furthermore resonant excitation modes.

\begin{figure}[!ht]

\vspace{10mm}
\centering
 \vspace*{-5mm} \rotatebox{0}{\resizebox{!}{10
cm}{%
\includegraphics [scale=0.25]{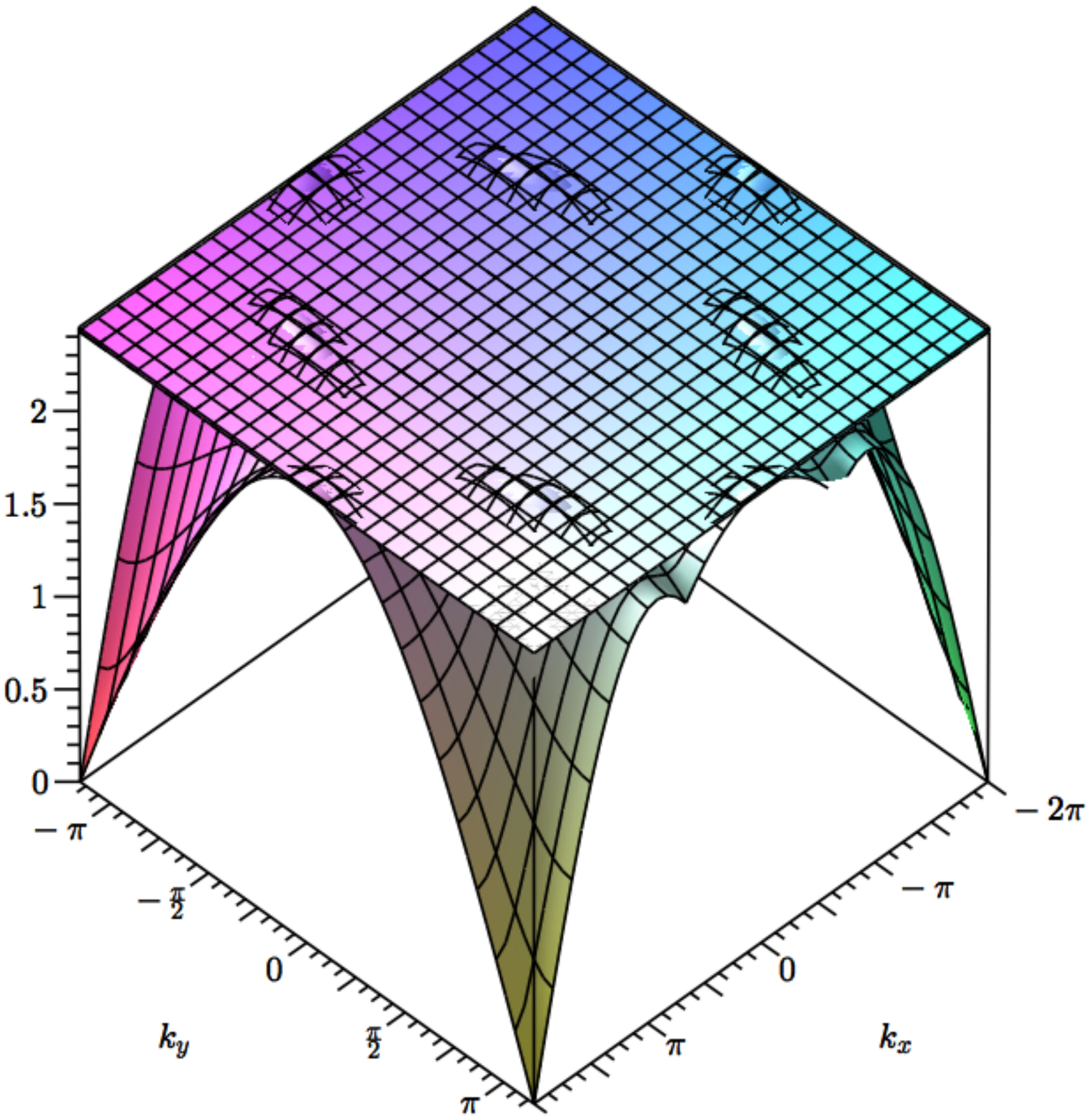} }} \\
\centerline{(a)} 
\vspace*{-5mm} \rotatebox{0}{\resizebox{!}{10
cm}{%
\includegraphics [scale=0.25]{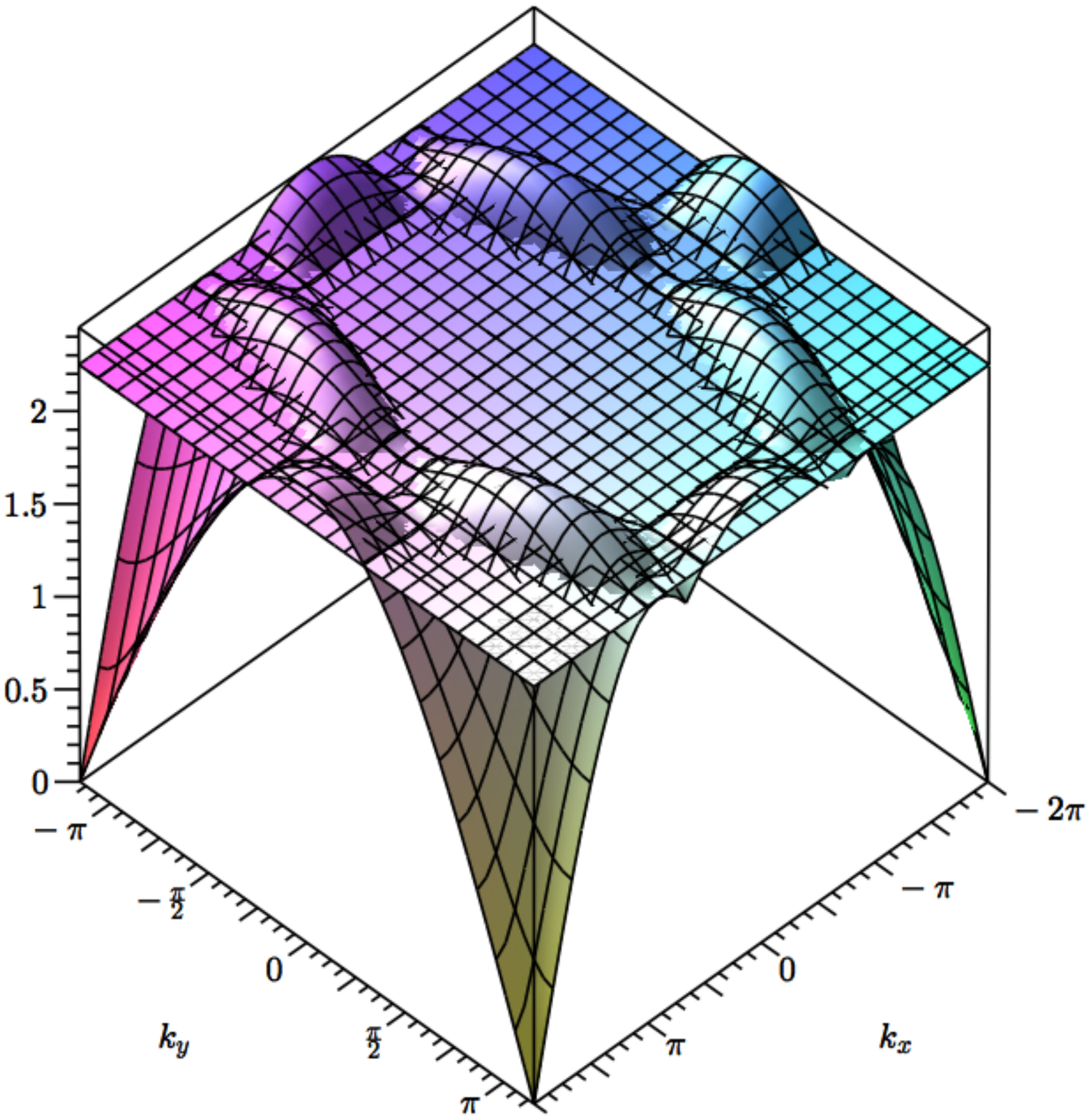}
\includegraphics [scale=0.25]{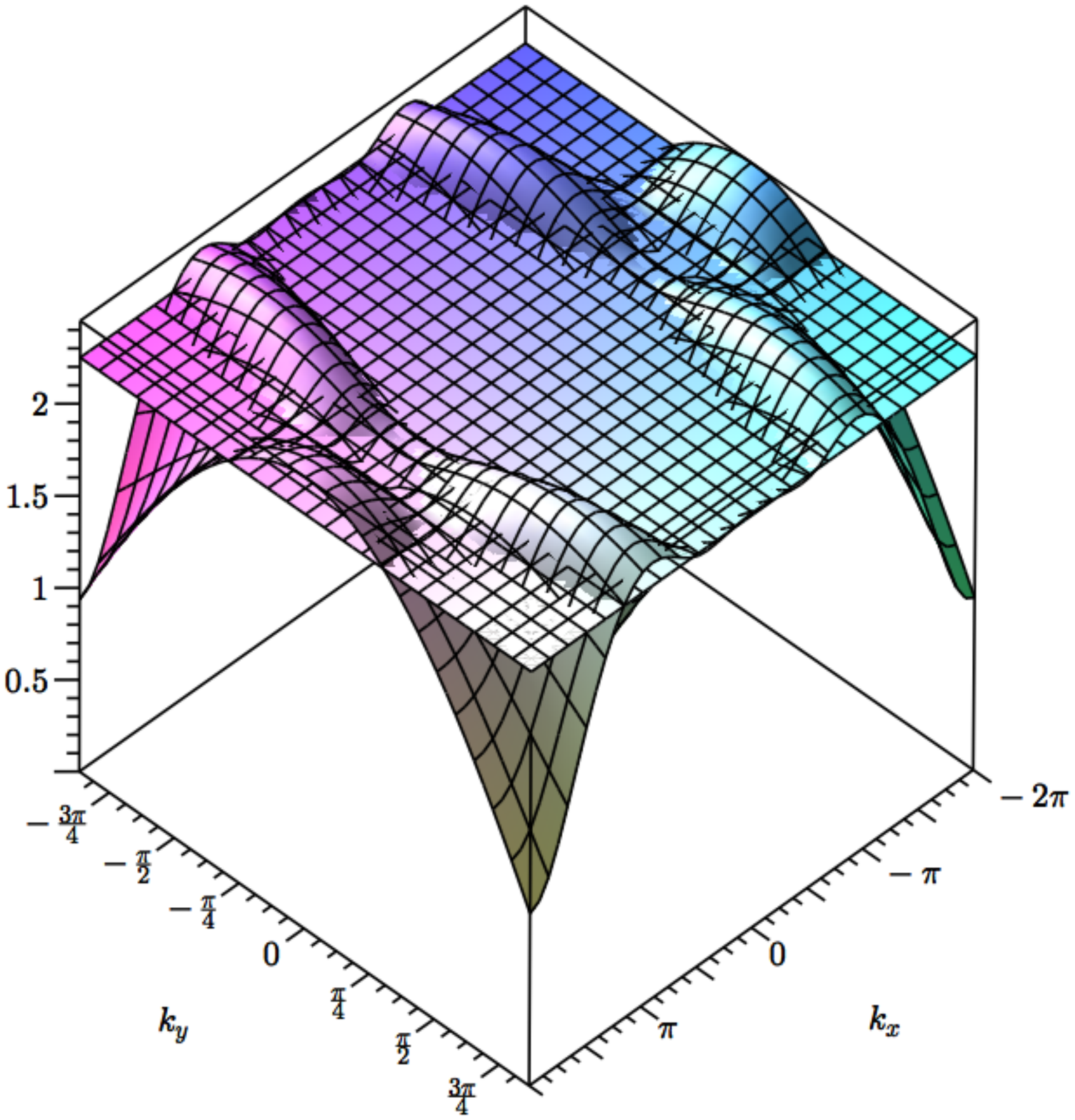}}} \\
\centerline{(b) ~~~~~~~~~~~~~~~~~~~~~~~~~~~~~~~~~~~~~~~~~~~~~~~~~~~~~~~~~(c)}

\vspace{5mm}
 \caption{(a) Dispersion surface $\Go=\Go_1(k_x,k_y)$  intersecting with the plane $\Go=\Go_1^*-0.02$. (b) Dispersion surface $\Go=\Go_1(k_x,k_y)$ intersecting with the plane $\Go=\Go_2^*$. (c) Dispersion surface $\Go=\Go_2(k_x,k_y)$ intersecting with the plane $\Go=\Go_3^*$.}
    \label{omega1-resonant-plane}
\end{figure}

 The dispersion surfaces crossed by
 fixed-frequency planes are presented in \fig{omega1-resonant-plane}a,b,c: 
 in the first of these diagrams, the plane is taken slightly below the points of maximum, $\Go=\Go_1^*-0.02$, to make the location of the points visible; in the diagrams (b) and (c) the critical points are the saddle points on the dispersion surfaces.

\clearpage


%
%
%
%

\begin{figure}[!ht]

\centering
\vspace*{10mm} \rotatebox{0}{\resizebox{!}{7
cm}{%
\includegraphics [scale=0.20]{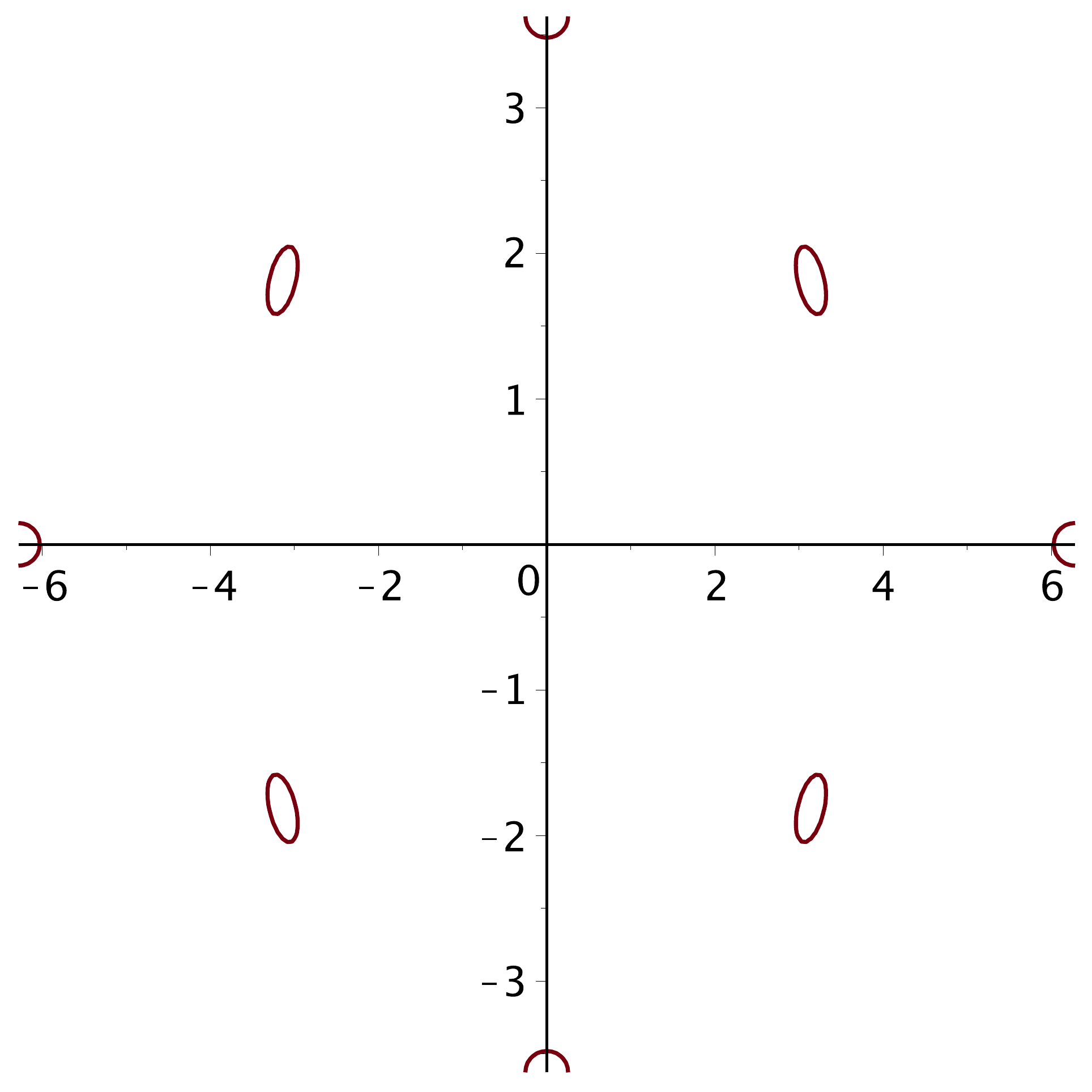}}}

\vspace{0mm}
 \caption{The trace of the $\Go_1$-surface crossed by the plane $\Go=\Go_1^*-0.02$.}
    \label{omega1-contour-1}
\end{figure}

\begin{figure}[!ht]

\centering
\vspace*{10mm} \rotatebox{0}{\resizebox{!}{7
cm}{%
\includegraphics [scale=0.20]{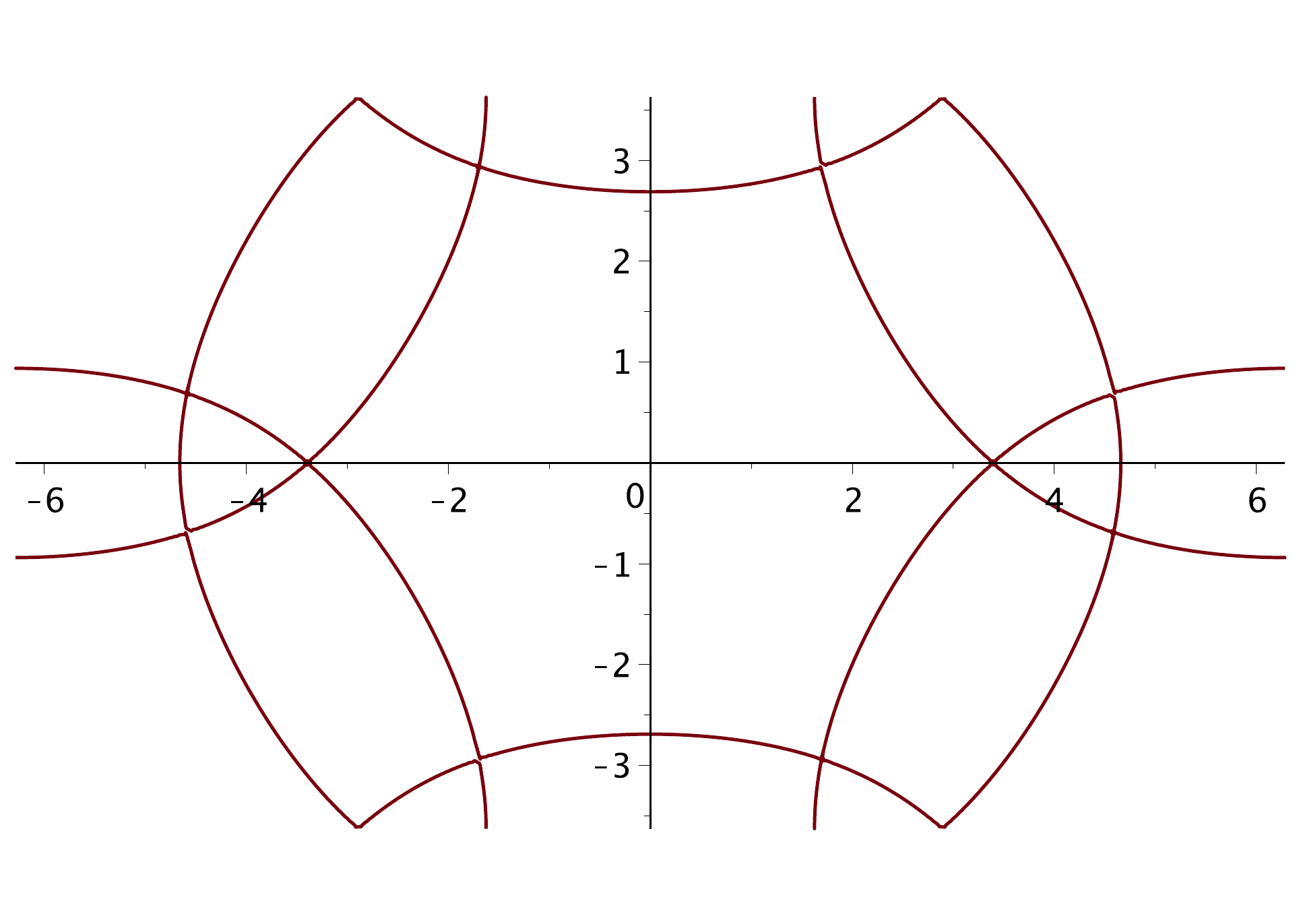}}}

\vspace{10mm}
 \caption{Resonant-frequency slowness diagram, $\Go_1=\Go_2^*=2.25$.}
    \label{omega1-contour-2}
\end{figure}

\begin{figure}[!ht]

\centering
\vspace*{10mm} \rotatebox{0}{\resizebox{!}{7
cm}{%
\includegraphics [scale=0.20]{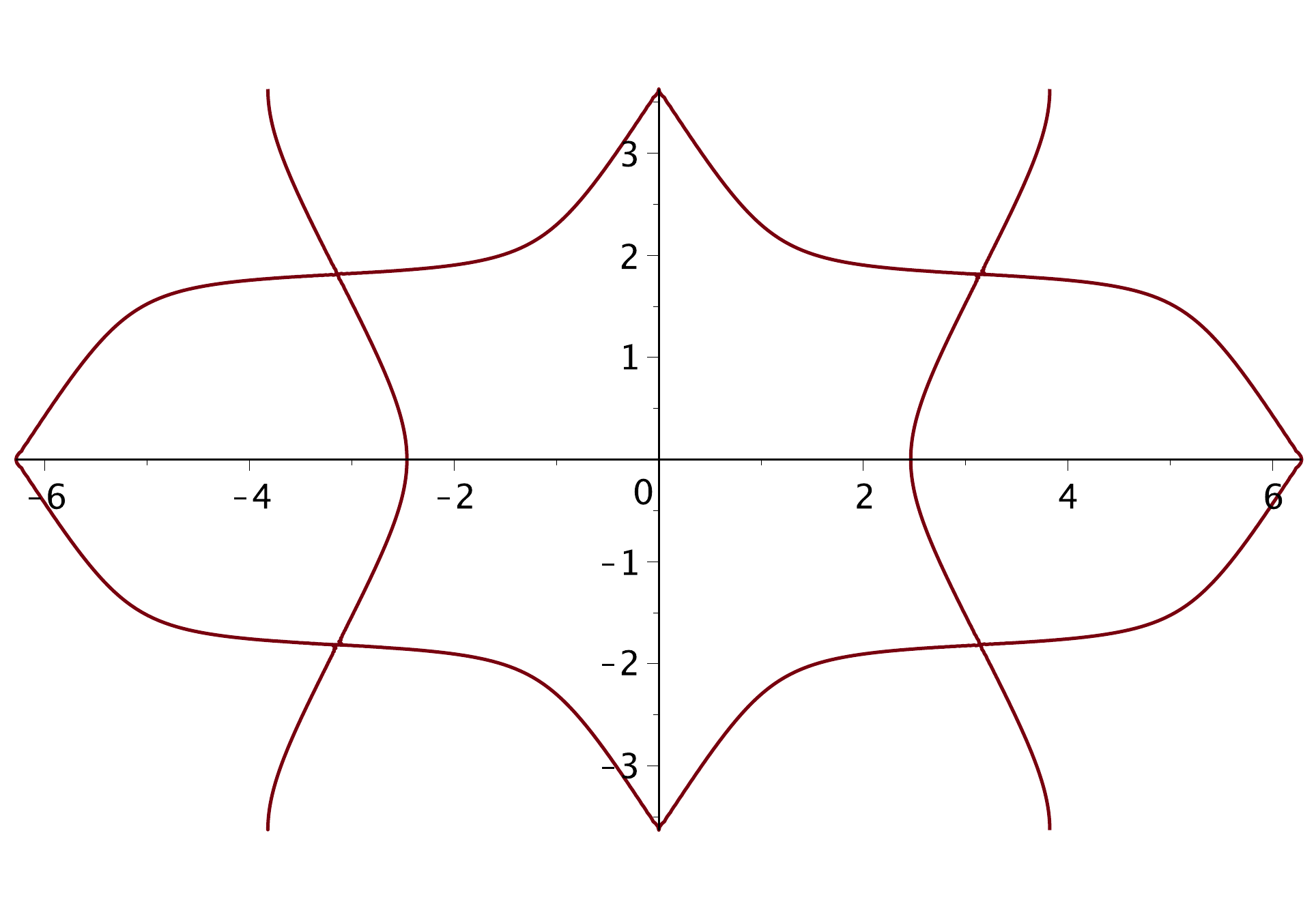}}}

\vspace{0mm}
 \caption{Resonant-frequency slowness diagram, $\Go_2=\Go_3^*=\sqrt{2}$.}
    \label{omega2-contour}
\end{figure}

The  resonant-frequency dispersion contours are plotted in \fig{omega1-contour-1} $-$ \fig{omega2-contour}. It is observed that the saddle points of the dispersion surfaces correspond to  the crossing points of the slowness contour diagram.
As shown in Table 1 (also illustrated by slowness contours), there exist 8, 14 and 8  resonant points at the resonant frequencies, $\Go^*_1, \Go^*_2$ and $\Go^*_3$, \res. In particular, the critical points corresponding to the frequencies $\Go^*_2$ and $\Go^*_3$ are the saddle points on the dispersion surfaces, whereas $\Go=\Go^*_1$ represents the point of maximum.

\subsection{Resonant versus non-resonant excitations}

First, we evaluate $\bfm{B}_{\Go,\bfm{k}}$, as in \eq{10},  at the resonant points. The required values are presented in Table 2 below.

%
Recall that the latter equality in \eq{5} evidences that the dispersion surfaces intersect only at the origin, $\Go_1=\Go_2=0$, and at the corner points of the periodicity rectangle, $k_x=\pm 2\pi, k_y=\pm 2\pi/\sqrt{3}$, where $\Go_1=\Go_2=0$ as well. At the critical points representing standing waves the dispersion surfaces do not intersect, i.e. $\Go_1\ne\Go_2$.
We also note that no growing wave is excited by an action corresponding to zero values of $B_{xx}, B_{yy}$ or $B_{x,y}$ presented in Table 2.

\newpage

The representation \eq{9} is valid for non-resonant excitations, where $\GD\ne 0$. There is no such a steady state in the resonant case, and a transient problem should be considered. For the transient problem we have to replace $\GO$ by $-\p^2/\p t^2$. In this way, for the following considerations the Laplace transform on $t$ is useful, and a formulation for the resonant harmonic excitation is obtained substituting
$ \GO =-s^2$ and $\bfm{P}\E^{\I\Go t} \to \bfm{P}/(s-\I\Go)$, where $s$ is the Laplace transform parameter. Thus, we have
 \beq \bfm{u}^{FFL} = \f{\bfm{B}\bfm{P}}{(s-\I\Go)\GD}\,.\eeq{9a}

\begin{table}[htb]
\centering
\begin{tabular}{@{\qquad}c@{\qquad}c@{\qquad}c@{\qquad}c@{\qquad}c@{\qquad}}
\toprule
\multicolumn{1}{c}{} \\
$ ij$ & $B_{xx}$ & $B_{yy}$ & $B_{xy}$\\
\midrule
\multicolumn{1}{c}{
}\\
$11/12$ & $0$ & $-4$ & $0$ \\
 \midrule
 $13$ & $-3$ & $-1$ & $-\sqrt{3}$\\
 \midrule
$21/22$ & $-1.6875$ & $0$ & $0$\\
 \midrule
$23/24$ & $-0.421875$ & $-1.265625 $ &$-0.730709$\\
 \midrule
$31/32$ & $4$ & $0$  & $ 0$ \\
 \midrule
$33$ & $1$ & $3$ & $-\sqrt{3}$ \\
\bottomrule
\end{tabular}
\caption{\label{Table2}
Coefficients $B_{xx}, B_{yy}, B_{xy}$ in the representation \eq{10} of the elastic compliance.}
\end{table}

%

\subsection{The dispersion relations in vicinities of the resonant points   and the characteristic lines    }
 The quadratic terms of the power expansion of $\GO \,(\GO_1$ or $\GO_2)$
 \beq \GO - (\Go_i^*)^2 \sim a_{ij}q_x^2 +b_{ij}q_y^2+c_{ij}q_x q_y~~~(k_{x,y}-k_{x,y;ij}=q_{x,y;ij}\to 0) \eeq{ktitpe1}
for neighbourhoods of the resonant points located in the first quadrant, $0\le k_{x;ij}\le 2\pi, 0\le k_{y;ij}\le 2\pi/\sqrt{3}$, are presented in Table 3:

\begin{table}[htb]
\centering
\begin{tabular}{@{\qquad}c@{\qquad}c@{\qquad}c@{\qquad}c@{\qquad}c@{\qquad}}
\toprule
\multicolumn{1}{c}{} \\
$ ij$ & $a_{ij} $ & $b_{ij}$ & $c_{ij}$\\
\midrule
\multicolumn{1}{c}{
}\\
$11/12$ &  $-0.375$ & $-1.125$ & $0$\\
 \midrule
$13$  & $-0.9375$ & $-0.5625$ & $-0.649519$\\
 \midrule
$21/22$ & $-0.984375$ & $1.265625  $ & $ 0 $\\
 \midrule
$23/24$ & $ 0.703125  $ & $-0.421875  $ & $ -1.94856 $\\
 \midrule
$31,32$ & $0.875 $ & $ -0.375  $ & $0 $ \\
 \midrule
$33$ & $ -0.0625  $ & $0.5625  $ & $-1.082532 $ \\
\bottomrule
\end{tabular}
\caption{\label{Table3}
Coefficients of the quadratic approximation \eq{ktitpe1} in the vicinity of critical points of dispersion surfaces.}
\end{table}



The polynomial in \eq{ktitpe1} corresponds to the differential operator
 \beq E=a_{ij}\f{\p^2}{\p x^2} +b_{ij}\f{\p^2}{\p y^2}+c_{ij}\f{\p^2}{\p x\p y}\,.\eeq{tpcte1}
The resonant points on the dispersion surfaces at $\Go= \Go_{2,3}^*$ are the saddle points, where the coefficients $a_{ij}$ and $b_{ij}$ are nonzero and different by the sign. Thus, the equation $E=0$ is hyperbolic at these frequencies. It is satisfied by an arbitrary function of $\Gf x-y$ with
 \beq \Gf= \f{c}{2a}\pm\sqrt{\gl(\f{c}{2a}\gr)^2 -\f{b}{a}}\,,~~~ (a, b)=( a_{ij},b_{ij})\,.\eeq{tpcte2}
In addition, if $c_{ij}\ne 0$, due to the symmetry there exist another saddle point with $c_{ij1}=-c_{ij}$. This results in additional  values of the parameter $\Gf$ with a different sign of the first term in \eq{tpcte2}. Therefore the slopes of
the characteristic lines  are given by the angles relative to the $x$-axis
 \beq \Ga = \pm \arctan\gl( \f{c}{2a}\pm\sqrt{\gl(\f{c}{2a}\gr)^2 -\f{b}{a}}\gr)\,.\eeq{phichls}
The values of such angles for the  rays associated with the saddle points of the dispersion surfaces are summarised  in Table 4.

\begin{table}[htb]
\centering
\begin{tabular}{@{\qquad}c@{\qquad}c@{\qquad}c@{\qquad}c@{\qquad}c@{\qquad}}
\toprule
\multicolumn{1}{c}{} \\
$ij$ & $\Ga$\\
\midrule
\multicolumn{1}{c}{
}\\
$21,22$ & $\pm (\pi/3-\Gl_2), \pm (2\pi/3 +\Gl_2)~~~(\Gl_2=0.19913547) $ \\
 \midrule
$23/24$ &  $\pm\Gl_2, \pi \pm\Gl_2, \pm (\pi/3+\Gl_2), \pm (2\pi/3 -\Gl_2) $\\
 \midrule
$31,32$ & $\pm (\pi/3-\Gl_3), \pm (2\pi/3 +\Gl_3)~~~(\Gl_3=0.46755781) $\\
 \midrule
$33$ & $\pm\Gl_3, \pi \pm\Gl_3, \pm (\pi/3+\Gl_3), \pm (2\pi/3 -\Gl_3) $  \\
\bottomrule
\end{tabular}
\caption{\label{Table4}
Coefficients of the quadratic approximation \eq{ktitpe1} in the vicinity of critical points of dispersion surfaces.}
\end{table}

%

It follows that the critical rays for $\Go=\Go_2^*$ and $\Go= \Go_3^*$ are oriented symmetrically with respect to each lattice bond line as should be, namely
 \beq \Ga = \pm 0.19913547 +\f{\pi n}{3}~(\Go=\Go_2^*)\,,~~\Ga = \pm 0.46755781 +\f{\pi n}{3}~(\Go=\Go_3^*)\,,~~n=0,1,...,5\,,\eeq{crda1}
where the directional angle $\Ga=0$ corresponds to a bond line.

\subsection{Green's tensor asymptotics}
Substituting $s=\I\Go+s'$ in \eq{9a} we obtain the following expressions for the oscillation amplitudes
 \beq \bfm{U}^{FFL}= \gl(\bfm{u}\E^{-\I\Go t}\gr)^{FFL} = \f{\bfm{B}\bfm{P}}{s'\GD}\eeq{9b}
 with
 \beq \GD = [(s'+\I\Go)^2+\GO_1][(s'+\I\Go)^2+\GO_2]\,.\eeq{10b}

To proceed we now note that for a resonant frequency, $\Go=\Go_i^*$, the amplitude growing in time is defined by integration in \eq{11} in the vicinities of the resonant points where $\GO_1$ or $\GO_2$ is equal to $(\Go_{i}^*)^2$.
We note that $\GO-(\Go_{i}^*)^2\to 0$ as $\bfm{k}$ tends to a resonant point. It follows from this that for the asymptotic representation of the resonant oscillation amplitude we need to consider the case as  $s'\to 0$. In this way, putting $(s'+\I\Go)^2 \sim -\Go^2 +2\I\Go s'$ we obtain the following asymptotic relation of $\GD$ associated with a resonant point
 \beq \GD\sim C[s' - \I(a q_x^2+b q_y^2+c q_x q_y)]\,,\eeq{12}
where the constant $C$ is different for different resonant frequencies, i.e.
 \beq C=2\I\Go(\GO_1-\Go^2)~~~(\Go=\Go_3^*)\,,~~~C=2\I\Go(\GO_2-\Go^2)~~~(\Go=\Go_{1,2}^*)\,.\eeq{c12j}
In turn, the coefficients $a, b, c$ are summarised in Table 3.

Using the relation in \eq{12} for the asymptotic regime, $t\to\infty, \bfm{x}/t\to 0$, the integration in the Fourier inverse transform can be extended over the whole $\bfm{k}$ plane. It follows that the wave associated with a saddle point $(ij)$, associated with the resonance frequency $\Go^*_i$, is asymptotically defined as
 \beq
 &&\bfm{u}_{ij} \sim \bfm{U}_{ij}\exp[\I(\Go_i^* t -\bfm{k}_{ij}  \cdot \bfm{x})]\,,\n
&& \mbox{and} \n
 &&\bfm{\dot{U}}_{ij}\sim\f{\sqrt{3}}{16\pi^2 C}V.p.\int_{-\infty}^\infty \int_{-\infty}^\infty \exp[\I (t(a q_x^2+bq_y^2+cq_xq_y)-xq_x-yq_y)]\D q_x\D q_y
 \bfm{B}_{ij}\bfm{P}\n
&& = \f{\sqrt{3}}{8\pi C\sqrt{c^2-4ab}\,\, t}\exp\gl(\I R\gr)\bfm{B}_{ij}\bfm{P}\,,~~~R=\f{bx^2+ay^2-cxy}{(c^2-4ab)t}\,,(a,b,c)=(a,b,c)_{ij}\,. \eeq{14}
It follows that the asymptotic approximation for the displacement amplitude is given by 
 \beq \bfm{U}_{ij} =\int_0^t \bfm{\dot{U}}_{ij} dt \sim  \f{\sqrt{3}}{8\pi C\sqrt{c^2-4ab}}\mbox{Ei}\gl(1, -\I R\gr)\bfm{B}_{ij}\bfm{P}\,,\eeq{Ei}
and
 \beq
 \mbox{Ei}(1,\I R) = - \mbox{Ci}(|R|)+\I\gl[\mbox{Si}(R)-\f{\pi}{2}\mbox{sign}R\gr]\n
 \sim -\ln |R| -\Gg -\I\f{\pi}{2}\mbox{sign}R~~~(|R|\ll 1)\,,\eeq{16}
where $\Gg$ is Euler's constant.

The above asymptotic solution is valid outside the characteristic line, $R=0$, where it becomes singular. To evaluate the wave amplitude on this line, one needs higher-order terms in the Taylor expansion of the dispersion surface near a critical point. Note that the characteristic lines defined by the equation $R=0$ coincide with those in  \eq{phichls}.

\begin{figure}[!ht]

\centering
\vspace*{-15mm} \rotatebox{0}{\resizebox{!}{9
cm}{%
\includegraphics [scale=0.25]{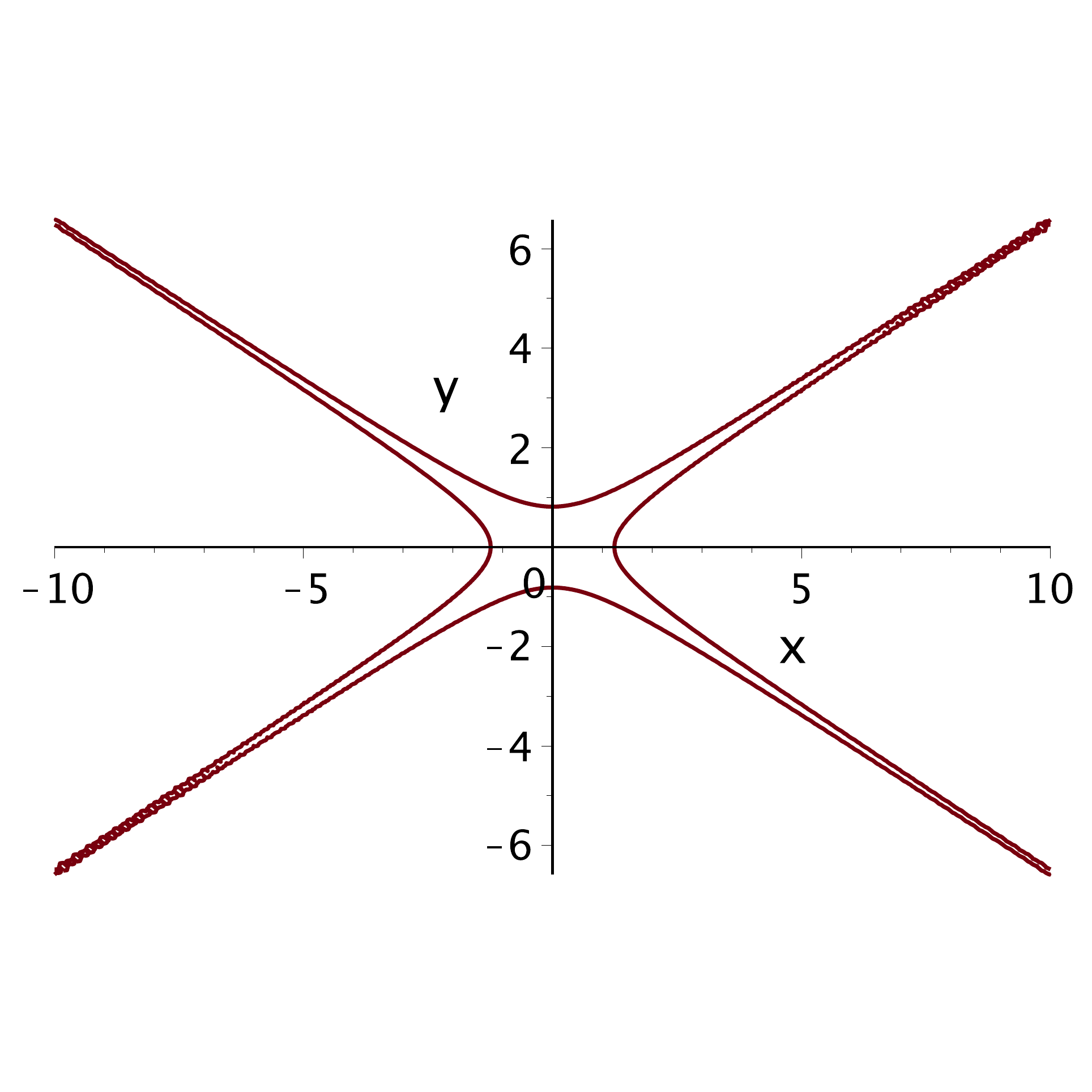}}}

\vspace{-10mm}
 \caption{The star-wave configuration associated with the resonant point $31$, $t=2000$, $|\mbox{Ei}(1,-\I R)|=8$.  }
    \label{star31}
\end{figure}

\begin{figure}[!ht]

\centering
\vspace*{-15mm} \rotatebox{0}{\resizebox{!}{10
cm}{%
\includegraphics [scale=0.25]{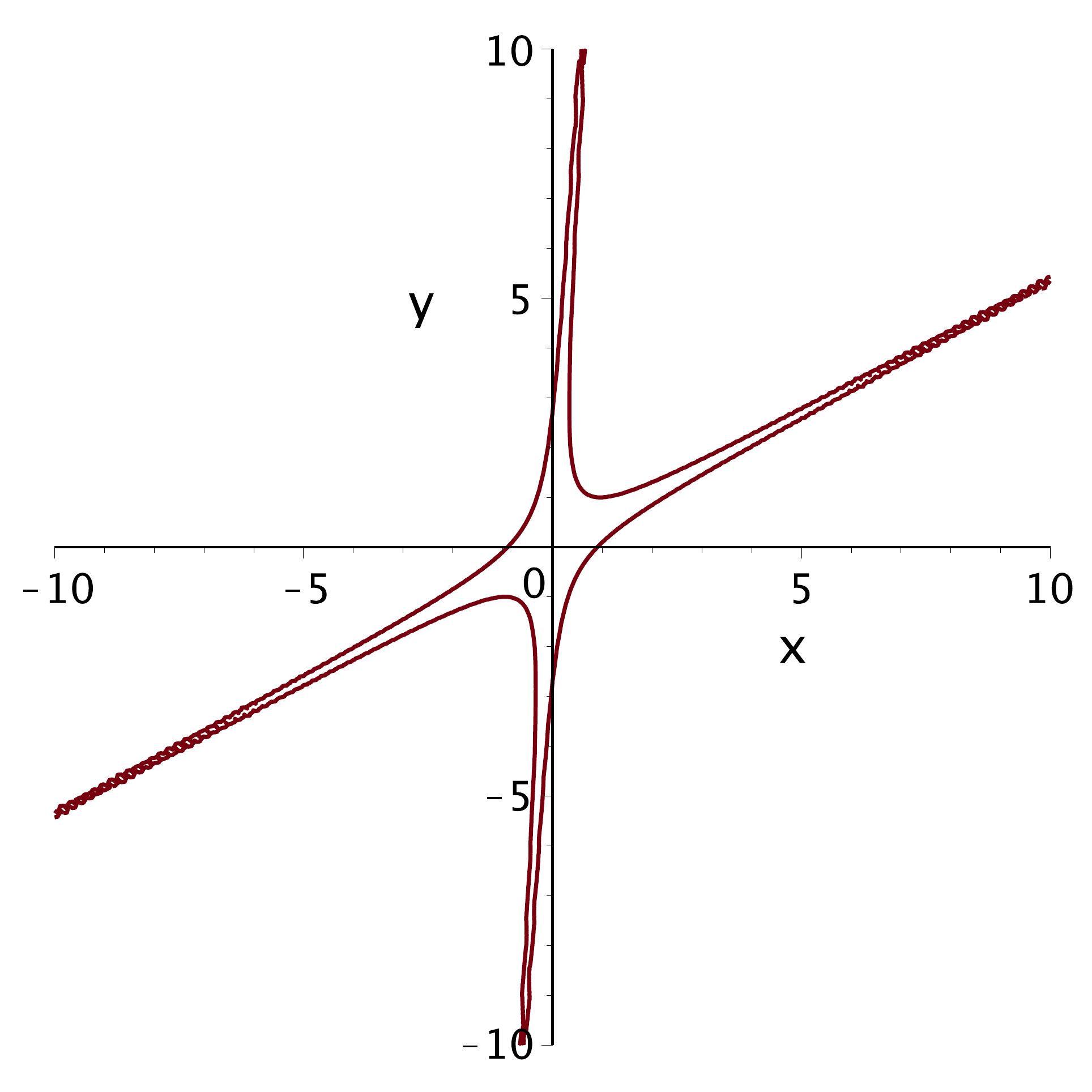}}}

\vspace{-5mm}
 \caption{The star-wave configuration associated with the resonant points $33$ in the first and third quadrants of the $k_x,k_y$-plane, $t=2000$, Ei(1,$-\I R)=8$. }
    \label{starC+33}
\end{figure}

\begin{figure}[!ht]

\centering
\vspace*{0mm} \rotatebox{0}{\resizebox{!}{10
cm}{%
\includegraphics [scale=0.25]{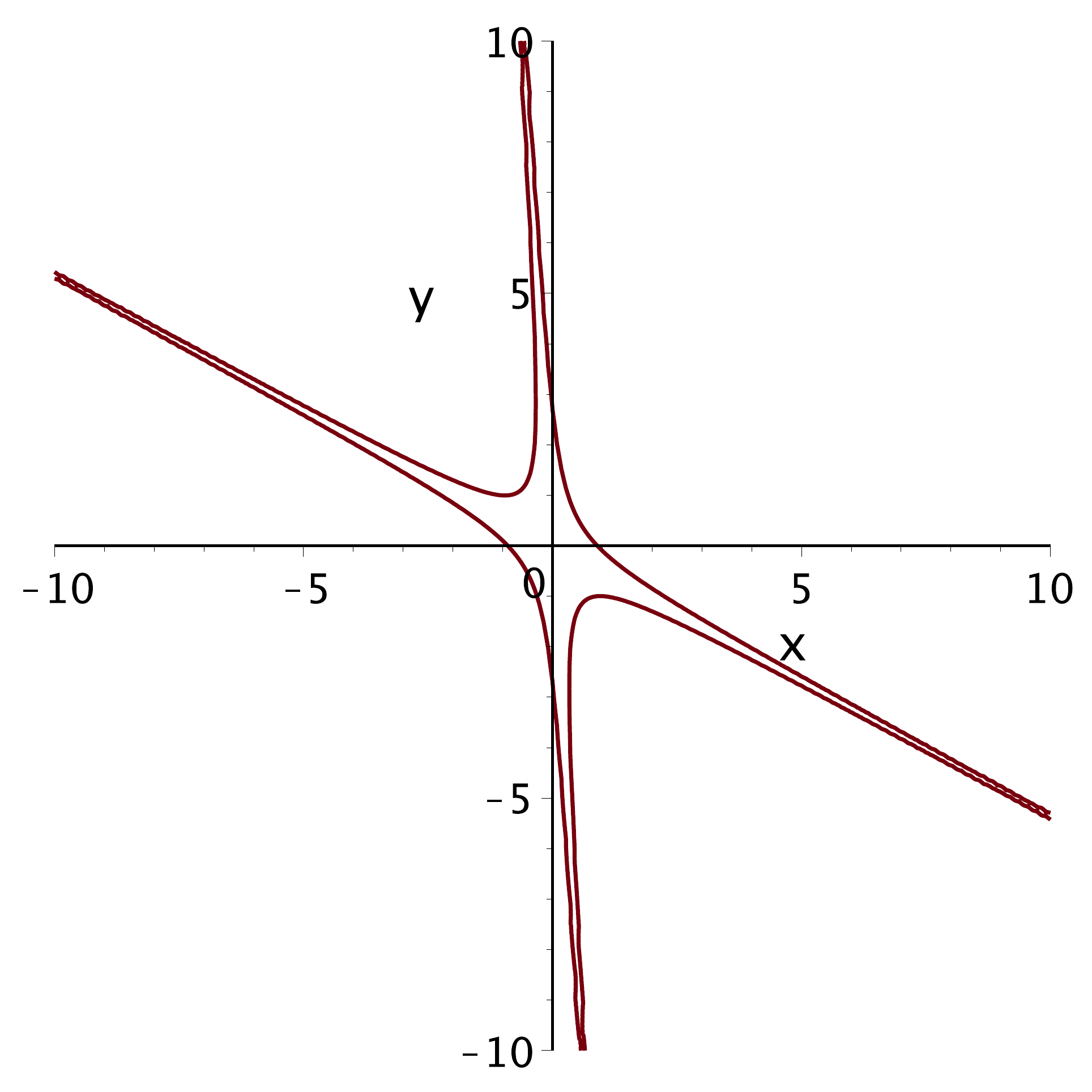}}}

\vspace{0mm}
 \caption{The star-wave configuration associated with the resonant point $\Go=\Go^*_{33}$ in the second and fourth quadrants of the $k_x,k_y$-plane, $t=2000$, Ei(1,$-\I R)=8$. }
    \label{starC-33}
\end{figure}

\begin{figure}[!ht]

\centering
\vspace*{-10mm} \rotatebox{0}{\resizebox{!}{9
cm}{%
\includegraphics [scale=0.20]{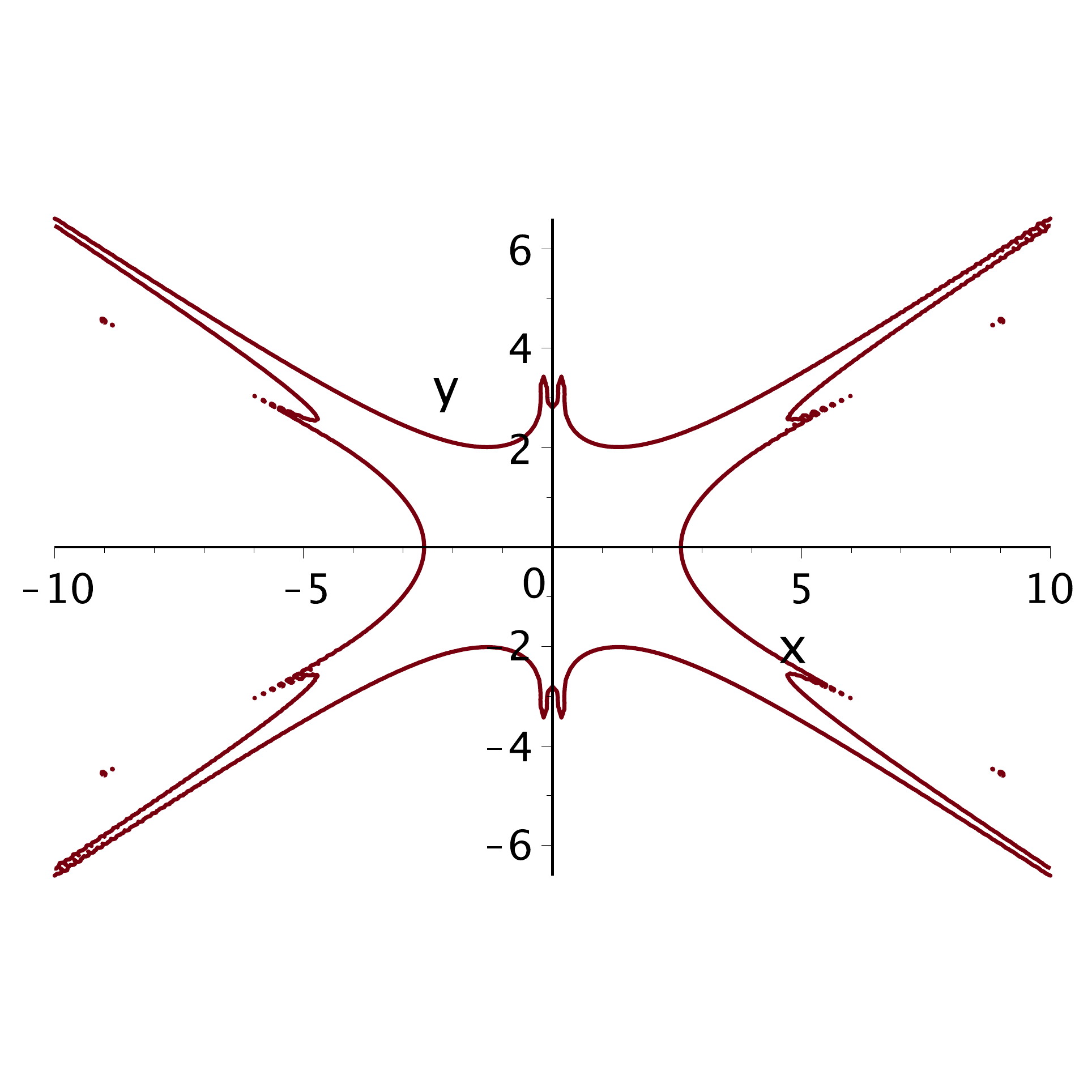}}}

\vspace{-10mm}
 \caption{The resonant wave excited by the horizontal unit force, $P_x=1, P_y=0$ at $\Go=\Go_3^*, t=2000, |\bfm{U}|=0.2.$}
    \label{horizontal}
\end{figure}

\begin{figure}[!ht]

\centering
\vspace*{10mm} \rotatebox{0}{\resizebox{!}{9
cm}{%
\includegraphics [scale=0.20]{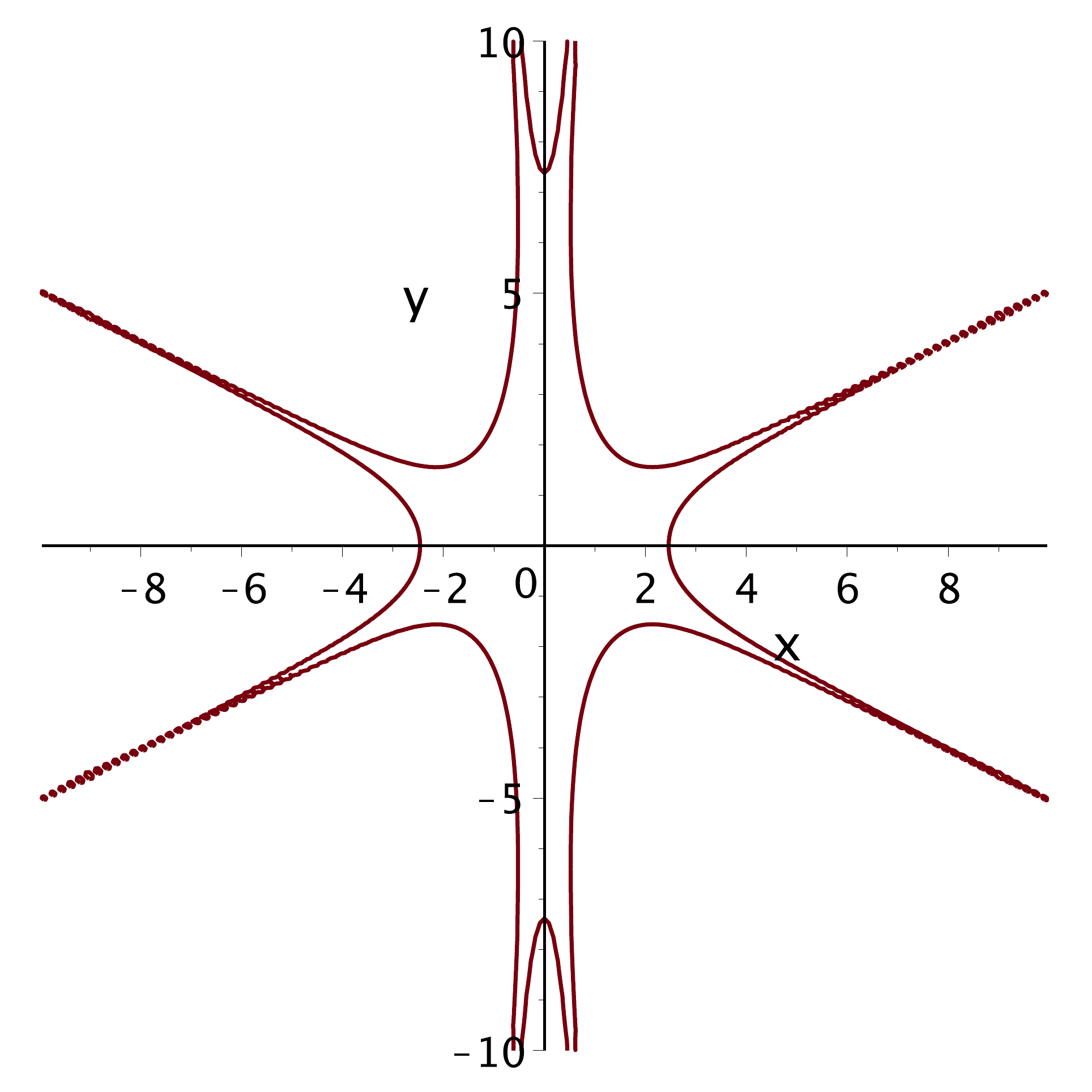}}}

\vspace{0mm}
 \caption{The resonant wave excited by the vertical unit force, $P_x=0, P_y=1$ at $\Go=\Go_3^*, t=2000, |\bfm{U}|=0.2.$}
    \label{vertical}
\end{figure}

A typical star-shaped wave form  associated with a saddle point is shown in \fig{star31}, where the level lines are plotted.
Note that for this symmetric wave
the coefficient $c_{ij}$ in  \eq{tpcte1} is equal to zero.
In \fig{starC+33}, \fig{starC-33}, we show the wave forms represented by the level lines of the amplitude of the displacement along the horizontal bond, $u_x$, for Ei$(1,R)=3$ at $t=400$. The amplitude is greater inside the star and lower outside it. The  resonant waves excited at $\Go=\Go_3^*$ are shown in \fig{horizontal} and \fig{vertical}.

\vspace{1in}

\section{Homogenisation and concluding remarks}
In accordance with \eq{tpcte1} and \eq{12}, the homogenized equation for $\bfm{\dot{U}}$ associated with a resonant point is
 \beq a\f{\p^2 \bfm{\dot{U}}}{\p x^2}+b \f{\p^2 \bfm{\dot{U}}}{\p y^2} +c\f{\p^2\bfm{\dot{U}}}{\p x\p y}-\I \f{\p\bfm{\dot{U}}}{\p t}=\f{\bfm{B}\bfm{P}}{C}\,, \eeq{18}
 where the coefficients, $a,b,c$, are generally different for different resonant points (see \eq{12}, \eq{c12j} and Table 3).
For a saddle resonant point, $ab<0$,  this equation is hyperbolic and the critical rays on the ($x,y$)-plane correspond to its characteristics. The solution to this equation is presented in \eq{14}, \eq{Ei}.

Recall that different homogenization corresponds to different resonant points at the same frequency. This does not allow for a global homogenisation, which could correspond to a given resonant frequency. Instead that initiates this hypothetic action. Instead, the global asymptotic Green's tensor is built, which represents the continuous approximation of the wave field.
In conclusion, note that while the Green's tensor is defined by the inverse transforms, the directions of the resonant wave localization can be seen in the resonant-frequency dispersion contours. Indeed, (a) the excited propagation waves correspond to the level lines since the other free waves have different frequencies; (b) the contribution to the resonant wave is not given by the resonant point itself but by a set of the waves corresponding to the level lines in a vicinity of the former; (c) the group velocities of the latter are directed along the normal to the level line. Thus, the `star' directions coincide with the normals to the level lines at the saddle resonant points, that is, at the cross- and angle points of the lines. Note that the energy transfer in 1D resonant waves is considered in Slepyan and Tsareva (1987).

\vspace{10mm}
\noindent {\bf Acknowledgements}

\noindent
The authors greatfully acknowledge support from the FP7 Marie Curie IAPP grant  \# 284544-PARM2.

\vskip 18pt
\begin{center}
{\bf  References}
\end{center}
\vskip 3pt

\inh Ayzenberg-Stepanenko, M., Slepyan, L.I., 2008. Resonant-Frequency Primitive Waveforms and Star Waves in Lattices.
Journal of Sound and Vibration 313, 812-821. \\ DOI: 10.1016/j.jsv.2007.11.047

\inh Bakhvalov, N.S., and Panasenko, G.P., 1984. Homogenization: Averaging Processes in Periodic Media. Mathematics and Its Applications (Soviet Series) 36. Kluwer Academic Publishers, Dordrecht-Boston-London.

\inh Bensoussan, A., Lions, J.L., and Papanicolaou, G.,  1978, Asymptotic Analysis for Periodic Structures. North-Holland, Amsterdam.

\inh Brun, M., Guenneau, S.R.L., Movchan, A.B., 2009, Achieving control of in-plane elastic waves, Appl. Phys. Lett., 94, 061903.

\inh Colquitt, D.J, Jones, I.S., Movchan, N.V., Movchan A.B., McPhedran, R.C., 2012, Dynamic anisotropy and localization in elastic lattice systems. Waves in Random and Complex Media, Vol 22, Issue 2, 143-159.

\inh Craster, R.V., Kaplunov, J., and Postnova, J., 2010. High-frequency asymptotics, homogenisation and localisation for lattices. Q.Jl Mech. Appl. Math. 63(4), 497-519.

\inh Craster, R.V., Nolde, E., and Rogerson, G.A., 2009. Mechanism for slow waves near cutoff frequencies in periodic waveguides. Phys. Rev. B79, 045129.

\inh Lord Rayleigh, On the influence of obstacles arranged in rectangular order upon the
properties of a medium, Phil. Mag. 34 (1892), pp. 481–502.

\inh Marchenko, V.A., and Khruslov, E.Y, 2006. Homogenization of Partial Differential Equations. Springer, Heidelberg.

\inh Milton, G.W., Briane M., Willis, J.R., 2006, On cloaking for elasticity and physical equations with a transformation invariant form. New Journal of Physics, 8, 248, doi: 10.1088/1367-2630/8/10/248.

\inh Milton, G.W. and Nicorovici, N.A., 2006, On the cloaking effects associated with anomalous localized resonance. Proc. Royal Soc., Ser. A, vol. 462, No 2074, 3027-3059.

\inh Norris, A.N, 2008, Acoustic cloaking theory. Proc. Royal Soc., Ser. A, vol. 464, No 2097, 2411-2434.

\inh Sanchez-Palencia, E., 1980, Non-Homogeneous Media and Vibration Theory, Springer Lecture Notes in Physics, Vol. 127, Heidelberg.

\inh Slepyan, L.I., and Tsareva, O.V., 1987. Energy Flux for Zero Group Velocity of the Carrying Wave. Sov. Phys. Dokl., 32, 522-524.

\inh 	Zhikov, V.V., Kozlov, S.M., Oleinik, O.A., 1994. Homogenization of differential operators and integral functionals,  Springer, Heidelberg.

\end{document}